\documentclass[journal,twoside,web]{ieeecolor}
\usepackage{generic}
\usepackage{cite}
\usepackage{amsmath,amssymb,amsfonts}
\usepackage{algorithmic}
\usepackage{comment,cuted,flushend}
\usepackage{graphicx,pgfplots}
\usepackage{textcomp}
\def\BibTeX{{\rm B\kern-.05em{\sc i\kern-.025em b}\kern-.08em
    T\kern-.1667em\lower.7ex\hbox{E}\kern-.125emX}}
\newcommand{\disp}{\displaystyle}
\newcommand{\z}{\boldsymbol{z}}
\newcommand{\x}{x}
\newcommand{\T}{\top}
\renewcommand{\L}{Q}
\renewcommand{\u}{\boldsymbol{u}}
\newcommand{\y}{\boldsymbol{y}}
\newcommand{\pd}[2]{\frac{\partial #1}{\partial #2}}
\renewcommand{\d}{{\rm d}}
\newcommand{\Ha}{\mathcal{H}}
\newcommand{\inner}[1]{\left\langle #1 \right\rangle}
\newcommand{\eig}{ {\rm eig}}

\newcommand{\R}{{\mathbb R}}

\newtheorem{theorem}{Theorem}
\newtheorem{lemma}{Lemma}

\begin{document}
\title{Exponential Decay Rate of Linear Port-Hamiltonian Systems. A Multiplier Approach}
\author{Luis A. Mora, \IEEEmembership{Member, IEEE}, and Kirsten Morris, \IEEEmembership{Fellow, IEEE}
\thanks{This research was supported by a Discovery Grant from NSERC and by a grant from the Faculty of Mathematics at the University of Waterloo.}
\thanks{Luis A. Mora and Kirsten Morris are with the Department of Applied Mathematics at University of Waterloo. 200 University Avenue West,	Waterloo, ON, Canada  N2L 3G1, (e-mails: lmora@uwaterloo.ca, kmorris@uwaterloo.ca).}
}

\maketitle

\begin{abstract}
	In this work, the multiplier method is extended  to obtain a general lower bound of the exponential decay rate in terms of the physical parameters for port-Hamiltonian systems in one space dimension with boundary dissipation. The physical parameters of the system may be spatially varying. It is shown that under assumptions of boundary or internal dissipation the system is exponentially stable. This is established through a Lyapunov function defined through a general multiplier function. Furthermore, an explicit bound on the decay rate in terms of the  physical parameters is obtained. The method is applied to a number of examples. 
\end{abstract}


\section{Introduction}\label{sec:introduction}
	\IEEEPARstart{E}xponential stability is a desirable property  of most  systems, including those modelled by partial differential equations. As shown in such works as \cite{Komornik1994}, \cite{Russell1994}, \cite{Xu1996} and \cite{Zuazua2012}  exponential decay of a partial differential equation through boundary control or dissipation is directly related to the exact observability of the system.  Furthermore, determining not only exponential stability but also an expression for  the exponential decay rate in terms of the system parameters is of theoretical interest, and also in practical performance such as analyzing control system performance.
An  important strategy to obtain an explicit expression for the exponential decay rate of dynamical infinite dimensional systems is the multiplier method \cite{Yan2022,Cheng2021,Rivera2021,Guesmia2003,Guo2004,Wu2014}. Several approaches of the multiplier method can be found in literature. For example, in \cite{Komornik1994} the system dynamics are multiplied by $m(\x)$ and the state variables, and  integrated in the space and time to derive the exponential decay rate of the state variable norm. Alternatively, the multiplier function is used to build an auxiliary Lyapunov functional whose exponential decay is related to the decay of the system energy; see the exposition in \cite{Tucsnak2009}.
	
A Port-Hamiltonian formulation of boundary-controlled distributed parameter systems was initially introduced in \cite{VanderSchaft2002,LeGorrec2005} and extended to problems with internal  dissipation in \cite{Villegas2006}. Sufficient conditions for the well-posedness of   linear PHS in one-dimensional spatial domains were established in  \cite{Zwart2010}. 
	Exponential stability of one-dimensional boundary controlled port-Hamiltonian systems has been studied in a number of works, including \cite{Villegas2009,Jacob2012,Augner2014,Macchelli2016,Trostorff2022}. In these works,  sufficient conditions that guarantee exponential stability are obtained. An explicit lower bound for the exponential decay rate of the energy  a Timoshenko beam with boundary  and internal dissipation was obtained in  \cite{Mattioni2022} .  

	In this work, we extend the multiplier method to obtain a general lower bound of the exponential decay rate in terms of the physical parameters for general port-Hamiltonian systems in one space dimension with boundary dissipation. The physical parameters of the system may be spatially varying. A  formal description of the  port-Hamiltonian systems under study is provided in Section \ref{sec:PHS}. The main results are presented in Section 3. We show that under assumptions of boundary or internal dissipation the system is exponentially stable. This is established by  considering general multiplier functions $m(\x). $   In previous work  \cite{MoraMorrisMTNS} a  linear multiplier function, $m(\x)=\x-a$, was used to obtain the an explicit bound on the  decay rate for port-Hamiltonian systems with  constant coefficients and  $P_0=0$, $G_0=0.$ This result was illustrated by obtaining  an explicit bound on the exponential decay rate of a boundary damped piezo-electric beam with magnetic effects.    The approach was extended in \cite{MoraMorrisCDC} to systems with $P_0\neq0$ and/or $G_0\neq0$, such as in a Timoshenko beam. However, the result in \cite{MoraMorrisCDC}  is restricted  to systems whose physical parameters satisfy several conditions. Here, by considering more general multiplier functions it is shown that a wider class of systems, including those with variable material parameters, are exponentially stable. Furthermore, an explicit bound on the decay rate in terms of the  physical parameters is obtained.  In Section 4, we apply the method to a number of examples. A preliminary analysis of Example A, the boundary-damped  wave equation with a linear multiplier function appeared in \cite{MoraMorrisMTNS}; here we compare the use of a linear  and exponential multiplier. The use of the linear multiplier function for the simple wave equation is well-known, here we not only compare different multiplier functions, but also, regarding the boundary damping as a control variable, the dissipation is chosen to optimize the decay rate. Example B applies our result to a wave equation with variable cross-section and material parameters.  In Example C, the Timoshenko beam, the result in \cite{MoraMorrisCDC} is extended to beams with general parameters. One lesson from these  examples is that the decay rate obtained depends on the choice of multiplier function.  Conclusions are presented on Section 5.

\section{Port-Hamiltonian systems (PHS)}\label{sec:PHS}

	Consider an one-dimensional spatial domain  $\Omega=\{ \x\in [a,b]\ \} \subset\mathbb{R} .$ Denote by $\z(\x,t)$ the $n$ state variables of a system on $\Omega$.   In this work, the following class of linear boundary controlled port-Hamiltonian systems \cite{Zwart2010} is considered:
	\begin{align}
		\pd{\z(\x,t)}{t}&-\left[P_1\pd{ }{\x}+\left[P_0-G_0\right]\right]\L(\x)\z(\x,t)=0\label{eq:PHS_1}
			\end{align}
	where $P_1=P_1^{\T}\in\mathbb{R}^{n\times n}$ is invertible, $\L(\x)=\L^{\T}(\x)>0\in H^1([a,b],\mathbb{R}^{n\times n})$, $P_0=-P_0^{\T}\in\mathbb{R}^{n\times n}$, $G_0=G_0^{\T}\geq0\in\mathbb{R}^{n\times n}$, and  $\disp\frac{n}{2}\times n$ real matrices $W_1$, $W_2$ and $\tilde{W}_1.$ Defining
	\begin{align}
	\u_b(t)&=W_1 \L(b)\z(b,t), \label{eq:u:b} \\
		\y_b(t)&=\tilde{W}_1 \L(b)\z(b,t)\label{eq:y_b} \\
		\u_a(t) &=W_2\L(a)\z(a,t) \label{eq:BC_a} ,
		\end{align}
		the boundary conditions are, for  some  $K=K^{\T}>0\in\mathbb{R}^{\frac{n}{2}\times \frac{n}{2}}$
			\begin{equation}
		\u_a(t)=0,  \quad
		\u_b(t)+K\y_b(t)=0 \label{eq:BC_1}
	\end{equation}
	or equivalently, 
	$$W_2\L(a)\z(a,t) =0,$$
	$$ W_1 \L(b)\z(b,t)+ K \tilde{W}_1 \L(b)\z(b,t)= 0 .$$
	The dissipative  boundary condition at $x=b$ can arise through natural boundary dissipation \cite[e.g.]{ZLMV} or as a controlled feedback with a measurement  $\y_b $ and controlled input $\u_b .$ 
	
	It will be assumed throughout that  $W_1$, $W_2$ and $\tilde{W}_1$ satisfy the following rank conditions.
	\begin{align}
		{\rm rank} \left(\begin{bmatrix}
			0 & W_2\\
			W_1 & 0
		\end{bmatrix}\right)=&n \quad\text{and}\label{eq:rank1}\\
		{\rm rank}\left(\begin{bmatrix}
			W_1\\\tilde{W}_1
		\end{bmatrix}\right)=n .
		\label{eq:ranks}
	\end{align}
	This guarantees that system \eqref{eq:PHS_1} defines a well-posed control system  \cite[Theorem 2.4]{Zwart2010}.

	The total energy of system \eqref{eq:PHS_1}-\eqref{eq:BC_1}  is 
	\begin{align}
		\Ha(t)=&\int_{a}^{b}\frac{1}{2}\z^\T(\x,t)\L(\x)\z(\x,t)\d\x\label{eq:H} .
	\end{align}
	Since $\L (\x )>0$ for all $x,$ this defines a norm on $L^2 ([a,b],\R^n)$ equivalent to the standard norm. 
	Differentiating \eqref{eq:H} along system trajectories,  assuming that $W_1^\T\tilde{W}_1=-W_2^\T\tilde{W}_2=P_1$ for some $\tilde{W}_2\in\mathbb{R}^{\frac{n}{2}\times n}$ and defining $\eta_K=\min \eig(K),$
	\begin{align}
		\frac{\d \Ha(t)}{\d t}=&-\int_{a}^{b}\z^\T(\x,t)\L G_0 \L\z(\x,t)\d\x \nonumber \\ & \quad +\frac{1}{2}	(\u_b^T (t)  \y_b (t)  + \y_b^T (t)   \u_b  (t)   ) 
		\nonumber\\
		\leq&-c_1\Ha(t)-\eta_K \|  \y_b  (t) \|^2 
		\label{eq:dtH_1}
	\end{align}
	where $c_1>0$ if internal dissipation $G_0 >0$  and $c_1=0$ otherwise. Exponential stability of the system when  $G_0 $ is not positive definite is not obvious. 
	
	\section{Exponential stability}

In this section the multiplier approach, see for example \cite{Tucsnak2009}, is modified and applied to the class of systems described in the previous section to obtain a explicit expression for the exponential decay rate in terms of the system parameters. 

\begin{lemma}\label{lemma:ExpSta}
	Let $\z(\x,t)\in L^2([a,b],\mathbb{R}^n)$ be the state of system \eqref{eq:PHS_1} on interval $\x\in[a,b]\subset\mathbb{R}$ and $\Ha(t)$ be the corresponding energy functional. If there exists a scalar functional $w(t)$ on $[a,b]$ of the state vector $\z(\x,t)$ such that
	\begin{align}
		|w(t)|\leq \frac{1}{\varepsilon_0}\Ha(t)
	\end{align}
	and
	\begin{align}
		\frac{\d w(t)}{\d t}\leq -\frac{1}{\varepsilon_1}\frac{\d \Ha(t)}{\d t}-c \Ha(t)
		\label{eq-wp}
	\end{align}
	for some  positive $\varepsilon_0$, $\varepsilon_1$ and $c, $ then $\Ha(t)$ decays exponentially. 
	Furthermore, defining $\disp M= \frac{\varepsilon_0+\varepsilon}{\varepsilon_0-\varepsilon}$ and decay rate $\disp \alpha=\frac{c\varepsilon\varepsilon_0}{\varepsilon_0+\varepsilon}$, $\forall \varepsilon\in[0,\min\{\varepsilon_0,\varepsilon_1\}],$
	\begin{align*}
		\Ha(t)\leq M e^{-\alpha t}\Ha(0) .
	\end{align*}
\end{lemma}

\begin{proof}
	Define $V_{\varepsilon}(t)=\Ha(t)+\varepsilon w(t)$, where $\varepsilon\in\mathbb{R}.$ 
	Note
	$$\Ha(t)-\varepsilon|w(t)|\leq V_{\varepsilon}(t)\leq \Ha(t)+\varepsilon|w(t)|.$$
	Also, since  $\displaystyle |w(t)|\leq \frac{1}{\varepsilon_0}\Ha(t),$
	\begin{align}
		\left(1-\frac{\varepsilon}{\varepsilon_0}\right) \Ha(t)\leq V_{\varepsilon}(t) \leq \left(1+\frac{\varepsilon}{\varepsilon_0}\right)\Ha(t)\label{eq:ExpDec0}
	\end{align}
	guaranteeng that $V_{\varepsilon}(t)$ is non-negative for all $\varepsilon\in[0,\varepsilon_0]$.
	
	Furthermore, using \eqref{eq-wp},  
	\begin{align*}
		\frac{\d V_{\varepsilon}(t)}{\d t}=&\frac{\d \Ha(t)}{\d t}+\varepsilon\frac{\d w(t)}{\d t}\\
		\leq&\left(1-\frac{\varepsilon}{\varepsilon_1}\right)\frac{\d \Ha(t)}{\d t} -\varepsilon c \Ha(t)
	\end{align*}
	  For any $\varepsilon\leq \varepsilon_1 ,$ \eqref{eq:ExpDec0} implies 
	\begin{align*}
		\frac{\d V_{\varepsilon}(t)}{\d t}\leq& -\varepsilon c \Ha(t) \leq -\frac{c\varepsilon}{1+\varepsilon/\varepsilon_0}V_{\varepsilon}(t)
	\end{align*}
	obtaining that $V_{\varepsilon}(t)=V_{\varepsilon}(0)e^{-\alpha t}$ where $\alpha=\dfrac{c\varepsilon\varepsilon_0}{\varepsilon_0+\varepsilon}$. Using again \eqref{eq:ExpDec0}, we obtain that $V_{\varepsilon}(0)\leq \left(1+\frac{\varepsilon}{\varepsilon_0}\right)\Ha(0)$ and $\Ha(t)\leq \frac{1}{1-\varepsilon/\varepsilon_0}V_{\varepsilon}(t)$. As a consequence,
	\begin{align}
		\Ha(t)\leq \frac{\varepsilon_0+\varepsilon}{\varepsilon_0-\varepsilon}e^{-\alpha t}
	\end{align}
	for all $\varepsilon\in[0,\min(\varepsilon_0,\varepsilon_1)]$, completing the proof.
\end{proof}


\begin{table}
	\centering
	\caption{System parameters}
	\label{tab:params}
	\begin{tabular}{cl}
		\hline Parameter  & Description\\ \hline
		$\disp\mu_\L$&$\disp\max_{\x\in[a,b]} \max \eig(\L(\x))$ \\
		$\disp\mu_B$&$\disp\max_{\x\in[a,b]} \max \eig(B(\x))$ \\
		$\disp\mu_{\Psi}$&$\disp\max_{\x\in[a,b]}\max \eig(\Psi(\x))$ \\
		$\disp\mu_{P_1}$&$\disp\sqrt{\max \eig(P_1^{-2})}$ \\
		$\disp \mu_m$&$\disp\max_{\x\in[a,b]} m(\x)$ \\
		$\disp\eta_\L$&$\disp\min_{\x\in[a,b]} \eig(\L(\x))$ \\
		$\disp\eta_K$&$\disp\min \eig(K)$\\
		\hline \multicolumn{2}{l}{{\footnotesize Matrices $B(\x)$ and $\Psi(\x)$ are defined in \eqref{eq:B} and \eqref{eq:Psi1}, respectively}.}
	\end{tabular}
\end{table}

\begin{theorem}\label{thm:1}
	Consider the port-Hamiltonian system with boundary dissipation given by \eqref{eq:PHS_1}-\eqref{eq:BC_1}. Define 
	\begin{align}
		\Psi(\x)=&\begin{bmatrix}
			-K\\ I
		\end{bmatrix}^\T \begin{bmatrix}
			W_1\\\tilde{W}_1
		\end{bmatrix}^{-\T}\L^{-1}(\x)\begin{bmatrix}
			W_1\\\tilde{W}_1
		\end{bmatrix}^{-1}\begin{bmatrix}
			-K\\ I
		\end{bmatrix} , \label{eq:Psi1}\\
		B(\x)=&\pd{\L(\x)}{\x}-\L(\x)(P_0+G_0)P_1^{-1}\nonumber\\&+P_1^{-1}(P_0-G_0)\L(\x)\label{eq:B} ,\\
		A_s(\x)=&¸\pd{m(\x)}{\x}\L(\x)-m(\x)B(\x)\label{eq:As1} .
	\end{align}
	Also for some $m(\x ) \in C([a,b])$ define the auxiliary function  of the state $\z$
	\begin{align}
		w(t)=\frac{1}{2}\int_{a}^{b} m(\x)\z^\T(\x,t)P_1^{-1}\z(\x,t)\d\x\label{eq:wt1}
	\end{align}
	Defining
	 $\disp\varepsilon_0=\dfrac{\eta_{\L}}{\mu_{m}\mu_{P_1}}$ and  $\disp\varepsilon_1=\frac{2\eta_{K}}{\mu_{m}\mu_{\Psi}},$
	if
	\begin{align}
		A_s(\x)>0\label{eq:As>0}
	\end{align}
	then
	for all $\varepsilon\in[0,\min\{\varepsilon_0,\varepsilon_1\}]$, 
		\begin{equation}\Ha (t) \leq M e^{-\alpha t }, \quad M = \frac{\varepsilon_0+\varepsilon}{\varepsilon_0-\varepsilon}  , \quad  \alpha=\frac{c\varepsilon\varepsilon_0}{\varepsilon+\varepsilon_0} \, . \label{H-bound}
		\end{equation}

\end{theorem}
\begin{proof}
	Using the Cauchy-Schwarz inequality,
	\begin{align*}
		|w(t)|=&\frac{1}{2}\left|\inner{m(\x)\z(\x,t),P_1^{-1}\z(\x,t)}\right|\\
		\leq& \frac{1}{2}\|m(\x)\z(\x,t)\|_{L^2}\|P_1^{-1}\z(\x,t)\|_{L^2}\\
		\leq&\frac{\mu_{m}\mu_{P_1}}{2}\|\z(\x,t)\|_{L^2}^2\leq\frac{\mu_{m}\mu_{P_1}}{\eta_{\L}}\Ha(t)
	\end{align*}
	Thus, $|w(t)|\leq \dfrac{1}{\varepsilon_0}\Ha(t).$	
	Similarly,
	\begin{align}
		\frac{\d w(t)}{\d t}=&\int_{a}^{b}m(\x)\z^\T(\x,t)P_1^{-1}\pd{\z(\x,t)}{t} ~\d\x \nonumber \\
		=&\int_{a}^{b}m(\x)\z^\T(\x,t)P_1^{-1}(P_0-G_0)\L(\x)\z(\x,t) ~\d\x\nonumber\\&+ \int_{a}^{b}m(\x)\z^\T(\x,t)\pd{\L(\x)\z(\x,t)}{\x} ~\d\x  \label{eq-dw}
	\end{align}
	Using the identity
	\begin{align*}
		\frac{1}{2}\pd{}{\x}\left(m(\x)\z(\x,t)^T \L(x)\z(\x)\right)=m(\x)\z^\T(\x,t)\pd{\L(\x)\z(\x,t)}{\x}\\+\frac{1}{2}\z^{\T}(\x,t)\left(\pd{m(\x)}{\x}\L-m(\x)\pd{\L(\x)}{\x}\right)\z(\x,t)
	\end{align*}
	\eqref{eq-dw}  is rewritten as
	\begin{align*}
		\frac{\d w(t)}{\d t}
		=&\left.\frac{1}{2}m(\x)\z^\T(\x,t)\L(\x)\z(\x,t)\right|_{a}^b\\&-\frac{1}{2}\int_{a}^{b}\z^\T(\x,t)A_s(\x)\z(\x,t)~\d\x
	\end{align*}
	where $A_s$ is defined in  \eqref{eq:As1}. Since $A_s$ is assumed positive, there exists a $c>0$ such that
	$$A_s(\x)\geq c\L(\x)>0.$$
	This implies that  
	\begin{align*}
		\frac{\d w(t)}{\d t}
		\leq&\frac{1}{2}m(b)\z^\T(b,t)\L(b)\z(b,t)\\&-\frac{c}{2}\int_{a}^{b}\z^\T(\x,t)\L(\x)\z(\x,t)~\d\x . 
	\end{align*}
	
	By assumption  \eqref{eq:ranks} $\begin{bmatrix}
		W_1\\\tilde{W}_1
	\end{bmatrix}$ is full rank and so $\disp \L(b)\z(b,t)=\begin{bmatrix}
		W_1\\\tilde{W}_1
	\end{bmatrix}^{-1}\begin{bmatrix}
		\u_b(t)\\\y_b(t)
	\end{bmatrix}$. Then, including the boundary dissipation \eqref{eq:BC_1} leads to, recalling  the definition of $\Psi$ in \eqref{eq:Psi1},
	\begin{align*}
		\frac{\d w(t)}{\d t}
		\leq&\frac{1}{2}m(b)\y^\T_b(t)\Psi\y_b(t)-c\Ha(t)\\
		\leq& \frac{1}{2}\mu_m\mu_{\Psi}|\y_b|^2-c\Ha(t)\\
		\leq& -\frac{1}{\varepsilon_1}\frac{\d \Ha(t)}{\d t}-c\Ha(t) .
	\end{align*}	 
	where $\varepsilon_1=\dfrac{2\eta_K}{\mu_{m}\mu_{\Psi}}.$ Lemma \ref{lemma:ExpSta} then implies the  bound on the exponential decay of  $\Ha(t)$ in  \eqref{H-bound}.
\end{proof}

\begin{lemma}\label{lemma:1}
	Consider the matrices $A_s(\x)$ and $B(\x)$ defined in \eqref{eq:As1} and \eqref{eq:B}, respectively. 
	%
	Defining $\disp m(\x)=Ce^{\beta (\x-a)}, $  if $\beta$ is sufficiently large then 
	\begin{align}
		A_s>0 ~,\forall \x\in[a,b]\label{eq:As>0_1} .
	\end{align}
\end{lemma}

\begin{proof}
	 With $m(\x)=Ce^{\beta (\x-a)}$, the matrix $A_s$ can be rewritten as
	\begin{align*}
		A_s=m(\x)\left(\beta \L(\x)-B(\x) \right) .
	\end{align*}
	
	Since $m(\x)>0$,  condition \eqref{eq:As>0_1} is satisfied if  matrix $\beta \L(\x)-B(\x)$ is  positive; that is  if $\disp \inf_{\x\in[a,b]} \eig\left(\beta \L(\x)-B(\x)\right)>0. $  
	Since $\disp\eta_{\L}=\inf_{\x\in[a,b]} \eig\left(\L(\x)\right) >0$ and recalling $\disp\mu_B=\sup_{\x\in[a,b]} \eig\left(B(\x)\right),$ if $\beta$ is chosen large enough that 
		\begin{align*}
		(\beta\eta_{\L}-\mu_B)>0 
	\end{align*}
then the required condition is satisfied. 
\end{proof}

Exponential stability of the class of systems described in section 2 now follows immediately, along with a bound on the decay rate. 
Lemma \ref{lemma:1} implies that exists at least one multiplier function, $m(\x)$, such that condition \eqref{eq:As>0_1} holds and so the system  \eqref{eq:PHS_1},\eqref{eq:BC_1} is exponentially stable.  Furthermore, Theorem \ref{thm:1} can be used to obtain a lower bound of the exponential decay rate for all systems with the form \eqref{eq:PHS_1}-\eqref{eq:BC_1} on interval $\x\in[a,b]$. 

From Lemma \ref{lemma:1}, there are  definitions for $M(\varepsilon)$ and $\alpha(\varepsilon)$ for all $\varepsilon$ on the interval $[0,\min\{\varepsilon_0,\varepsilon_1\}]$, and Theorem \ref{thm:1} provides explicit expressions of $\varepsilon_0$ and $\varepsilon_1$ for system \eqref{eq:PHS_1}-\eqref{eq:BC_1}. 
Using  the parametrization $\varepsilon=\xi\min\{\varepsilon_0,\varepsilon_1\}$ with $0<\xi<1$ leads to
\begin{align}
	M=&\begin{cases}
		\dfrac{1+\xi}{1-\xi} & \text{if}~ \varepsilon_0\leq \varepsilon_1\\
		\dfrac{\eta_{\L}\mu_{\Psi}+2\xi\eta_{K}\mu_{P_1}}{\eta_{\L}\mu_{\Psi}-2\xi\eta_{K}\mu_{P_1}} & \text{otherwise}
	\end{cases}\label{eq:M}\\
	\alpha=&\begin{cases}
		\dfrac{\xi}{\xi+1}\dfrac{c\eta_{\L}}{\mu_{P_1}\mu_{m}}& \text{if}~\varepsilon_0\leq \varepsilon_1\\
		\dfrac{2 c \eta_{K}\eta_{\L} \xi}{\mu_{m}\left(\eta_{\L}\mu_{\Psi}+2\xi\eta_{K}\mu_{P_1}\right)} & \text{otherwise}
	\end{cases}
\end{align}

Since $c$ and $\mu_{m}$ are affected by the multiplier function, an appropriate choice of $m(\x)$ improves the exponential decay rate bound obtained through Theorem \ref{thm:1}.  Considering $\disp m(\x)=Ce^{\beta (\x-a)}$ with $C>0$, as in the proof of Lemma \ref{lemma:1}, we obtain $m(a)=C$, $\mu_{m}=m(b)=Ce^{\beta (b-a)}$  and
\begin{align*}
	A_s=&m(\x)\left(\beta \L(\x)-B(\x) \right)\\ \geq& m(a)(\beta\eta_{\L}-\mu_B) \\ \geq  &c\L
\end{align*}
where $\disp c=\frac{C(\beta\eta_{\L}-\mu_B)}{\mu_{\L}}$. Then, the exponential decay rate bound is given by 
\begin{align}
	\alpha=&\begin{cases}
		\dfrac{\xi}{\xi+1}\dfrac{\eta_{\L}e^{-\beta(b-a)}(\beta\eta_{\L}-\mu_B)}{\mu_{\L}\mu_{P_1}} &\text{if}~ \varepsilon_0\leq \varepsilon_1\\
		\dfrac{2\xi\eta_{K}\eta_{\L} e^{-\beta(b-a)} (\beta\eta_{\L}-\mu_B)}{\mu_{\L}\left(\mu_{P_1}\eta_{K}+2\xi\eta_{\L}\mu_{\Psi}\right)} & \text{otherwise}
	\end{cases}\label{eq:alpha_beta}
\end{align}

\begin{theorem}The system \eqref{eq:PHS_1}-\eqref{eq:BC_1}  is exponentially stable. Furthermore, the decay rate is at least
	\begin{align}
		\alpha=&\begin{cases}
			\dfrac{\xi}{\xi+1}\dfrac{\eta^2_{\L} e^{-\left(1+\frac{ \mu_B}{\eta_{\L}}(b-a)\right)}}{(b-a)\mu_{\L}\mu_{P_1}} &\text{if}~ \varepsilon_0\leq \varepsilon_1\\
			\dfrac{2\xi\eta_{K}\eta_{\L}^2 e^{-\left(1+\frac{\mu_B}{\eta_{\L}}(b-a)\right)}}{(b-a)\mu_{\L}\left(\mu_{P_1}\eta_{K}+2\xi\eta_{\L}\mu_{\Psi}\right)} & \text{otherwise}
		\end{cases}\label{eq:op_alpha}
	\end{align}
\end{theorem}

\begin{proof}
Using the exponential multiplier function from Lemma \ref{lemma:1}  along with Theorem \ref{thm:1} yields the conclusion that the system is exponentially stable, along with bounds on $M$ and $\alpha.$
	Since $M$ is independent of $\beta$, the optimal decay rate is obtained by choosing $\beta$ to maximize $\alpha$, that is 
	\begin{align}
		\beta_{op}=\arg\max_{\beta} \alpha=\frac{\eta_{\L}+(b-a)\mu_B}{(b-a)\eta_{\L}}=\frac{1}{b-a}+\frac{\mu_B}{\eta_{\L}}
	\end{align}
	Then, the optimal decay rate is obtaining substituting $\beta_{op}$ in \eqref{eq:alpha_beta}. 
\end{proof}

As shown Lemma \ref{lemma:1}, choosing $m(\x)$ as an exponential function, condition \eqref{eq:As>0} can always be satisfied. Depending on the system, other options for $m(\x)$ may be possible. 
For example, consider the case $P_0=G_0=0$, $\L(\x)=L\x+D>0$, $\forall\x\in[a,b]$, with $D$ and $L$ defined positive. Choosing $m(\x)=q\x+d$ where $\disp \frac{q}{d}\geq -\frac{1}{a},$
matrix $A_s$ becomes
\begin{align}
	A_s=&\pd{m(\x)}{\x}\L(\x)-m(\x)\pd{\L(x)}{\x}\nonumber\\=&q(L\x+D)-(q\x+d)L=qD-dL \end{align}
Then, \eqref{eq:As>0} is satisfied if $\disp \frac{q}{d}>\max\left\lbrace\frac{\mu_L}{\eta_D} ,-\frac{1}{a}\right\rbrace$.  
This point is illustrated in an example in the next section. 

\section{Examples}


\subsection{Wave equation with boundary dissipation}
Consider the  wave equation in an one-dimensional spatial domain
\begin{align}
	\pd{}{t}\left(\rho\pd{w(\x,t)}{t}\right)=&\pd{}{\x}\left(\tau \pd{w(\x,t)}{\x}\right) & \forall\x\in[a,b] \label{eq:ExWave}
\end{align}
with boundary conditions
\begin{align}
	\pd{w(a,t)}{t}=&0 & \forall t\geq0\label{eq:EXWaveBC1}\\
	\tau\pd{w(b,t)}{\x}+k\pd{w(b,t)}{t}=&0& \forall t\geq 0 \label{eq:EXWaveBC2}
\end{align}
and $w(\x,0)\in L^2([a,b],\mathbb{R})$, where the density and elasticity parameters, $\rho$ and $\tau$ respectively, are constant.

Defining $\disp z_1=\pd{w(\x,t)}{\x}$ and $\disp z_2=\rho\pd{w(\x,t)}{t}$, the wave equation \eqref{eq:ExWave} is expressed as the port-Hamiltonian system
\begin{align}
	\pd{\z(\x,t)}{t}=&P_1\pd{}{\x}\left(\L\z(\x,t)\right), & \forall \x\in[a,b]\label{eq:WEPHS}
\end{align}
where $\z(\x,t)=\begin{bmatrix}
	z_1(\x,t) & z_2(\x,t)
\end{bmatrix}^\T$, $P_1=\begin{bmatrix}
	0 & 1\\1&0
\end{bmatrix}$ and $\L=\begin{bmatrix}
	\tau & 0\\0 & {1}/{\rho}
\end{bmatrix}$. Similarly, choosing
$W_1=W_2=\begin{bmatrix}
	1 &0
\end{bmatrix}$ and $\tilde{W}_1=\begin{bmatrix}
	0 &1
\end{bmatrix}$, the boundary conditions \eqref{eq:EXWaveBC1}-\eqref{eq:EXWaveBC2} can be rewritten in the form  \eqref{eq:u:b}-\eqref{eq:BC_a}.

\begin{figure}
	\centering
	\includegraphics[width=\columnwidth]{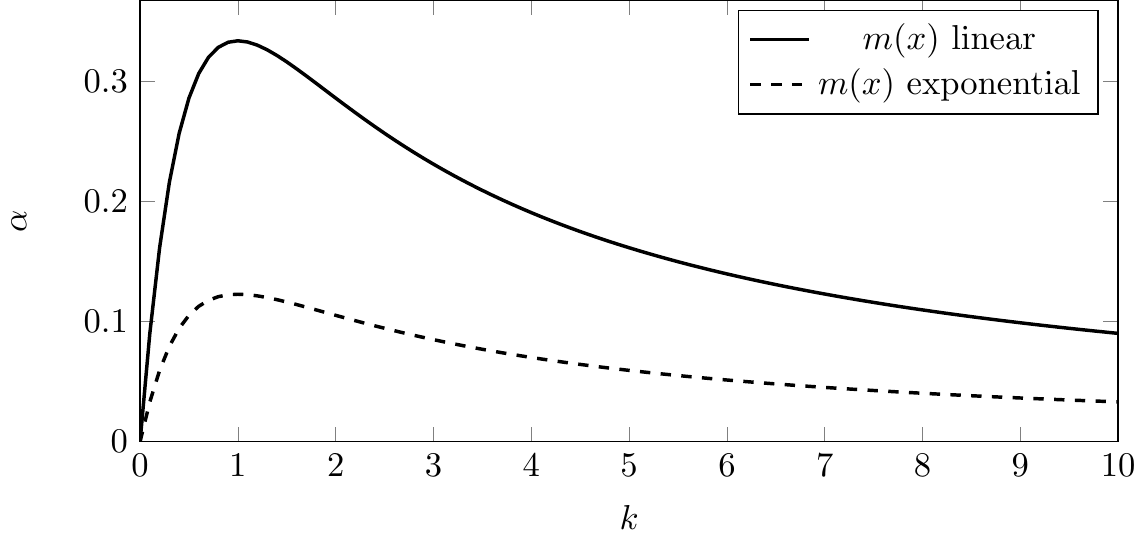}
	\caption{Bound on the exponential decay rate of the wave equation as a function of the boundary dissipation $k$ for different multiplier functions}
	\label{fig:WaveEq}
\end{figure}

Since $\L$ is a constant matrix and $P_0=G_0=0, $  $A_s=\pd{m(\x)}{\x}\L$. Condition \eqref{eq:As>0} holds if the multiplier function $m(\x)$ is monotonically increasing. Assuming unitary parameters, $\tau=\rho=1$, and spatial domain length, $b-a=1$, 
\begin{align*}
	\mu_B=&0&\mu_{P_1}=&\mu_{\L}=\eta_{\L}=1,&\mu_{\Psi}=&k^2+1 \\
	\varepsilon_0=&\frac{1}{\mu_m} & \varepsilon_1=&\frac{2k}{\mu_m (k^2+1)} & c=&\min_{\x\in[a,b]}\pd{m(\x)}{\x}
\end{align*}
Since $\disp \frac{2k}{k^2+1}\leq 1$ for all $k\geq 0, $ $\varepsilon_1\leq \varepsilon_0$ for any $k$.   Choosing $\disp\varepsilon=\frac{1}{2}\varepsilon_1$ we obtain $\disp M=\frac{k^2+k+1}{k^2-k+1}$ which is independent of the choice of $m(\x)$. 
With an exponential multiplier function, as in  Lemma \ref{lemma:1}, from \eqref{eq:op_alpha} we obtain that the decay rate $\alpha=\dfrac{ke^{-1}}{k^2+k+1}$. Alternatively, considering a linear multiplier function, $m(\x)=x-a$,    so $\mu_{m}=c=1$ and the exponential decay rate is $\alpha=\dfrac{k}{k^2+k+1}.$ This is   a better lower bound for the decay rate than the exponential multiplier function. This point is illustrated in  Figure \ref{fig:WaveEq}.

If the boundary dissipation comes from a control law, then the value of $k$ should be chosen to optimize the exponential decay rate. For this example the exponential decay rate $\disp \alpha=\frac{k}{k^2+k+1}$  is maximized  with $k=1. $ Note that this choice of $k$ is the same value for which no waves are reflected and the energy of the wave equation reaches zero in  finite time.

\subsection{Wave equation with  variable cross-section and boundary dissipation}
\begin{figure*}
	\centering
	\includegraphics[width=1.95\columnwidth]{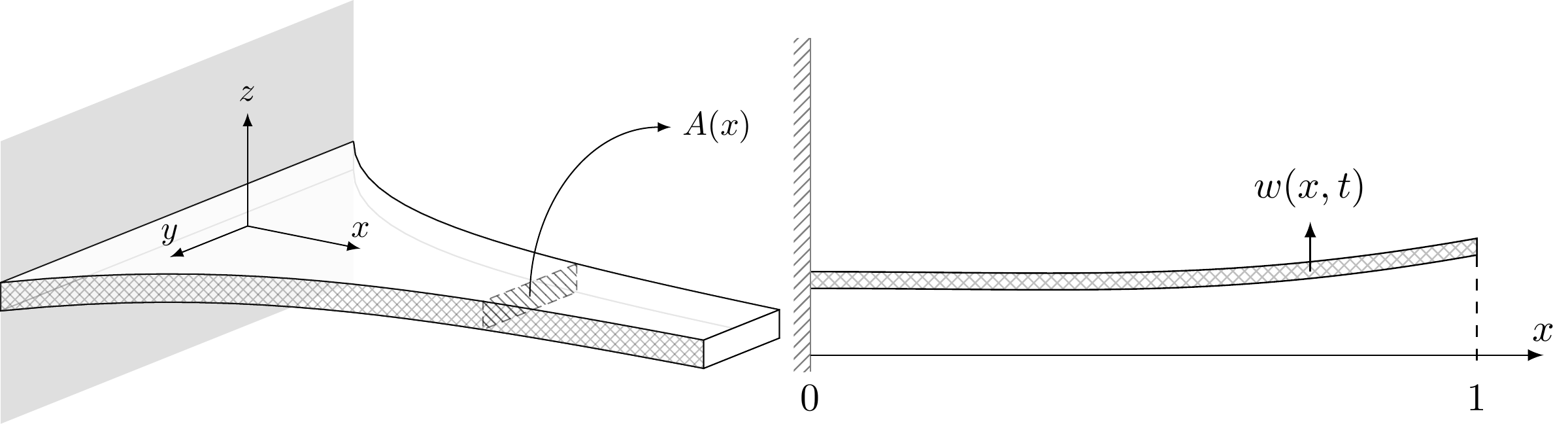}
	\caption{Vibrating string with non-uniform cross-sectional area}
	\label{fig:VS}
\end{figure*}

In this example, we consider a vibrating string with a non-uniform cross-sectional area $A(\x)$, as shown in Figure \ref{fig:VS}. The dynamics of the vertical displacements $w(\x,t)$ can be expressed as the wave equation with boundary dissipation described in \eqref{eq:ExWave}-\eqref{eq:EXWaveBC2}, where the physical parameters  $\rho(\x)=A(\x)\rho_0$ and $\tau(\x)=A(\x)\tau_0$ with $\rho_0$ and $\tau_0$ constant, and $w(\x,0)\in L^2([0,1],\mathbb{R})$. It will be assumed that $\tau_0=\rho_0=1$, boundary dissipation gain $k=0.5$ and
that cross-sectional area $$A(\x)=\dfrac{10-\x}{10}.$$
The port-Hamiltonian formulation of this vibrating string is similar to the previous analysis \eqref{eq:WEPHS} except that  
$$\L(\x)=\begin{bmatrix}
	\dfrac{(10-x)}{10} & 0\\
	0 & \dfrac{10}{(10-x)}
\end{bmatrix}.$$
As a consequence, 

	\begin{align}
		A_s(\x)=&\pd{m(\x)}{\x}\L(\x)-m(\x)\pd{\L(\x)}{\x}\nonumber\\
		=&\begin{bmatrix}
			\dfrac{m(\x)+\pd{m(\x)}{\x}(10-x)}{10} & 0\\
			0 & \dfrac{\pd{m(\x)}{\x}(10-x)-m(\x)}{0.1(10-x)^2}
		\end{bmatrix} . 
	\end{align}

Choosing the multiplier function $m(\x)=\x$, 
\begin{align*}
	A_s(\x)	=&\begin{bmatrix}
		1 & 0\\
		0 & \dfrac{10 \left(10-2\x\right)}{(10-x)^2}
	\end{bmatrix}>0, \quad\forall \x\in[0,1] \, . 
\end{align*}

The problem of finding the maximum $c$ such that $A_s(\x)\geq c\L(\x)$ is equivalent to finding the largest $c$ so that the eigenvalues of $A_s(\x)-c\L(\x)$ are non-negative; that is so the matrix
\begin{align*}
	\begin{bmatrix}
		1-c\dfrac{10-\x}{10} & 0\\
		0 & \dfrac{10}{(10-\x)}\left(\dfrac{10-2\x}{10-\x}-c\right)
	\end{bmatrix}
\end{align*}
is positive semi-definite. The largest such value of $c$ is $c=8/9$. 
Then,$$\disp\Psi(\x)=\frac{5}{2\left(10-\x\right)}+\frac{(10-\x)}{10},$$ 
\begin{align*}
	\mu_{\L}=&\frac{10}{9} & \eta_{\L}=&\frac{9}{10} & \mu_{\Psi}=&\frac{5}{4}\\
	\mu_{P_1}=&\mu_{m}=1 & 	\epsilon_0=&\frac{9}{10} & \varepsilon_1=&\frac{4}{5} \, . 
\end{align*}
Finally, choosing $\varepsilon=\varepsilon_1$ a bound on the exponential decay rate is 
$$\alpha=\frac{32}{85}\approx0.3765 . $$

\subsection{Timoshenko Beam}

Consider a Timoshenko beam  with variable material parameters on a bar $ x\in [a,b].$ Let   $\rho(\x)$, $\epsilon(\x)$ and $\iota(\x)$ be the mass per unit length,Young’s modulus, and moment of inertia of the cross section, respectively; $\iota_{\rho}(\x)=\iota(\x)\rho(\x)$ is the mass moment of inertia of the cross section; and  $\gamma(\x)$ and $\delta(\x)$ are the viscous damping coefficients. The shear modulus $\kappa(\x)= \xi G(\x)A(\x)$, where $G(\x)$ is the modulus of elasticity in shear, $A(\x)$ is the cross sectional area, and
$\xi$ is a constant depending on the shape of the cross section.  The parameters $k_1$ and $k_2$ are  boundary damping coefficients. This leads to the following partial differential equation
\begin{subequations}
	\begin{align}
		\rho(\x)\pd{^2w(\x,t)}{t^2}=&\pd{}{\x}\left(\kappa(\x)\left(\pd{w(\x,t)}{\x}-\phi(\x,t)\right)\right)\nonumber\\&-\gamma(\x)\pd{w(\x,t)}{t}\\
		\iota_{\rho}(\x)\pd{^2\phi(\x,t)}{t^2}=&\pd{}{\x}\left(\epsilon(\x)\iota(\x)\pd{\phi(\x,t)}{\x}\right)-\delta(\x)\pd{\phi(\x,t)}{t}\nonumber\\&+\kappa(\x)\left(\pd{w(\x,t)}{\x}-\phi(\x,t)\right)
	\end{align}\label{eq:TB_1}
\end{subequations} 
with boundary conditions
\begin{subequations}
	\begin{align}
		\pd{w(a,t)}{t}=\pd{\phi(a,t)}{t}=&0\\
		\kappa(b)\left(\pd{w(b,t)}{\x}-\phi(b,t)\right)+k_1\pd{w(b,t)}{t}=&0\\
		\epsilon(b)\iota(b)\pd{\phi(b,t)}{\x}+k_2\pd{\phi(b,t)}{t}=&0 \, . 
	\end{align}\label{eq:TB_2}
\end{subequations}
Set $z_1(\x,t)=\rho(\x)\pd{w(\x,t)}{t}$, $z_2(\x,t)=\iota_{\rho}(\x)\pd{\phi(\x,t)}{t}$, $z_3(\x,t)=\pd{w(\x,t)}{\x}-\phi(\x,t), $  $z_4(\x,t)=\pd{\phi(\x,t)}{\x}$ and $\z(\x,t)=\begin{bmatrix}
	z_1(\x,t) & z_2(\x,t) &z_3(\x,t)& z_4(\x,t)
\end{bmatrix}^\T . $   System \eqref{eq:TB_1} can be rewritten in the port-Hamiltonian formulation as in \cite{Mattioni2022} to obtain
$$\pd{\z(\x,t)}{t}-P_1\pd{\L(\x)\z(\x,t)}{\x}-\left[P_0-G_0\right]\L(\x)\z(\x,t)=0$$
where
\begin{align*}
	P_1=&\begin{bmatrix}
		0 & 0 & 1 & 0\\
		0 & 0 & 0 & 1\\
		1 & 0 & 0 & 0\\
		0 & 1 & 0 & 0
	\end{bmatrix}, \quad P_0=\begin{bmatrix}
		0 & 0 & 0 & 0\\
		0 & 0 & 1 & 0\\
		0 & -1 & 0 & 0\\
		0 & 0 & 0 & 0
	\end{bmatrix}, \\
	G_0=&\begin{bmatrix}
		\gamma & 0 &0 &0\\
		0 & \delta & 0 & 0\\
		0 & 0 & 0 & 0\\
		0 & 0 & 0 & 0
	\end{bmatrix}, \quad\text{and} \\
	\L(\x)=&\begin{bmatrix}
		\frac{1}{\rho(\x)} & 0 & 0 & 0\\
		0 & \frac{1}{\iota_{\rho}(\x)} & 0 & 0\\
		0 & 0 & \kappa(\x)& 0\\
		0 & 0 & 0 & \epsilon(\x)\iota(\x)
	\end{bmatrix} \, . 
\end{align*}

Similarly, defining
\begin{align}
	W_1=&\begin{bmatrix}
		0 & 0 & 1 & 0\\
		0 & 0 & 0 & 1\\
	\end{bmatrix}, \quad \tilde{W}_1=\begin{bmatrix}
		1 & 0 & 0 & 0\\
		0 & 1 & 0 & 0
	\end{bmatrix} \quad\text{and}\nonumber\\ W_2=&\begin{bmatrix}
		1 & 0 & 0 & 0\\
		0 & 1 & 0 & 0
	\end{bmatrix}
\end{align}
 the boundary conditions \eqref{eq:TB_2} can be written  in the standard form \eqref{eq:u:b}-\eqref{eq:BC_a} .

First, consider an inviscid beam with the same  parameters as in \cite{Mattioni2022}; that is $\gamma(\x)=\delta(\x)=0$, with  $\rho=0.2$kg/m, $\epsilon\iota=1.2\times 10^{-2}$Nm$^2$,  $\kappa=4\times10^{-3}$N, $\iota_{\rho}=2\times10^{-2}$kgm , $b-a=0.1$m and $k_1=k_2=k$. This leads to 
\begin{align*}
	\L=&\begin{bmatrix}
		5 & 0 & 0 & 0\\0 & 50 & 0 & 0\\0 & 0 & \dfrac{1}{250} & 0\\ 0 & 0 & 0 & \dfrac{1}{75}
	\end{bmatrix} \\ 
	B=&\begin{bmatrix}
		0 & -50 & 0 & 0\\-50 & 0& 0&0\\0 & 0 & 0 & \dfrac{1}{250}\\ 0 & 0 & \dfrac{1}{250}& 0
	\end{bmatrix} \\ 
	\Psi=&\begin{bmatrix}
		250k^2+\dfrac{1}{5} & 0\\0 & 75k^2+\dfrac{1}{50}
	\end{bmatrix}
\end{align*}
and $\mu_{P_1}=1$, $\mu_{\L}=\mu_B=50$, $\eta_{\L}=1/250$ and $\mu_{\Psi}=250k^2+1/5$.
Choosing a linear multiplier function, $m=x-a$ with $\mu_{m}=0.1$, then
\begin{align*}
	A_s(\x)=\begin{bmatrix}
		5 & 50(x-a)& 0 & 0\\50(x-a) & 50 & 0 & 0\\0 & 0 & \dfrac{1}{250}&\dfrac{a-x}{250}\\ 0 & 0 & \dfrac{a-x}{250} & \dfrac{1}{75}
	\end{bmatrix} 
\end{align*} 
which has eigenvalues 
\begin{align*}
	\eig(A_s(\x))=\begin{cases}
		\dfrac{13\pm\sqrt{36(x-a)^2+49}}{1500}\\ \\
		\dfrac{55\pm5\sqrt{400(x-a)^2+81}}{2}
	\end{cases} \, . 
\end{align*}
Since $0\leq x-a\leq 0.1$,  $\disp \min \eig(A_s(\x))\geq \frac{65-\sqrt{1234}}{7500}>3.9 \times10^{-3}$. 
This implies that     for sufficiently small $c>0, $  $\z^\T(\x,t) A_s(\x)\z(\x,t)\geq c\z^\T(\x,t) \L\z(\x,t)$. 
More precisely, the eigenvalues
	\begin{align*}
		\eig(A_s-c\L)=\begin{cases}
			\dfrac{\left(13 \pm 7 \sqrt{\left(\frac{6(x-a)}{7(1-c)}\right)^2+1}\right)(1-c)}{1500}
			\\ \\
			\dfrac{\left(55 \pm 45 \sqrt{\left(\frac{20(x-a)}{9(1-c)}\right)^2+1}\right)(1-c)}{2}
		\end{cases}
	\end{align*}
	need to be non-negatives. It is easy to check, through some simple calculations, that this condition is satisfied when $c\leq 1+\sqrt{\dfrac{1}{10}}$.
	
Thus, applying Theorem \ref{thm:1} with	$\varepsilon_0=\dfrac{1}{25}$, $ \varepsilon_1=\dfrac{100k}{1250k^2+1}$ and $c=0.6837$ and choosing $\varepsilon=\dfrac{1}{50}$ we obtain that $M=3$ and $\alpha=4.5 \times 10^{-3}.$ 
The bound on the decay rate in this example  is not improved with an exponential multiplier function.

Now we consider normalized physical parameters as in \cite{Mattioni2022}. That is, $\rho(\x)=\iota_{\rho}(\x)=\epsilon(\x)\iota(\x)=\kappa(\x)=\gamma=\delta=1$, boundary dissipation coefficients, $k_1=k_2=1$, and beam length $b-a=1.$  We obtain that $\eta_{K}=\eta_{\L}=\mu_{\L}=\mu_{P_1}=1$, $\mu_B=\sqrt{2}$ and $\mu_{\Psi}=2$. Then, considering a linear multiplier function, 
\begin{align*}
	A_s (\x ) =\begin{bmatrix}
		1&\x-a& \x-a&0\\
		\x-a&1&0 & \x-a\\
		\x-a&0 &1& a-\x\\
		0 & \x-a&a-\x&1
	\end{bmatrix}
\end{align*}
whose eigenvalues are $1\pm\sqrt{2}(\x-a)$. As a consequence, $A_s (\x)  >0$  only for $\x<a+\dfrac{1}{\sqrt{2}}<b$ and not in the entire interval $[a,b]$. The linear multiplier function $m(\x)=\x-a$, used with the previous set of parameters, cannot be used and it is necessary to consider another function.

Choosing an exponential multiplier function, as in Lemma \ref{lemma:1}, 
\begin{align*}
	\varepsilon_0=&\varepsilon_1=\frac{e^{-(1+\sqrt{2})}}{C}
\end{align*}
Then, varying $\xi$ on \eqref{eq:M} and \eqref{eq:op_alpha}, we obtain the values of $M$ and $\alpha$ shown in Table \ref{tab:1}. 

In  \cite{Mattioni2022} a Lyapunov approach is used for the stability analysis in the port-Hamiltonian formulation of a Timoshenko beam with viscous dissipation and unitary parameters, leading to an exponential decay rate of $0.0285$ with a $M=2.783.$
	Choosing $\xi=0.4713$, from \eqref{eq:op_alpha} we also obtain $M=2.783$ and $\alpha=0.0286.$
	
\begin{table}
	\caption{Values of $\alpha$ and $M$ for different choices of $\xi$}
	\label{tab:1}
	\centering
	\begin{tabular}{ccc}
		\hline $\xi$ & $M$ & $\alpha$\\\hline
		1/3 & 2 & 0.0224\\
		0.4713& 2.783 & 0.0286\\
		1/2 & 3 & 0.0298\\
		3/5 & 4 & 0.0335\\
		2/3 & 5 & 0.0358\\\hline
	\end{tabular}
\end{table}

%

\section{Conclusions}

An explicit formulation in terms of physical parameters for the exponential energy decay lower bound of a class of port-Hamiltonian systems with boundary dissipation on one-dimensional spatial domains have been presented. The choice of an exponential function, $m(\x)=Ce^{\beta(\x-a)}$, leads to a conclusion that provided that the boundary dissipation $K>0$  the system is exponentially stable. Furthermore, a
a lower bound on the decay rate is obtained $\alpha . $ This result applies to systems with variable physical parameters, as illustrated by several examples. 

For uniform systems,  $m(\x)$ is commonly chosen as a linear function; that is $m(\x)=\x-\x_0$, where $\x_0$ is chosen  so that $m(a)\geq 0$; see for example, \cite{Komornik1994,Tucsnak2009}.
 This choice of $m(\x)$ also works for uniform  port-Hamiltonian systems \eqref{eq:PHS_1}-\eqref{eq:BC_1} with $P_0=G_0=0$, as was shown in   \cite{MoraMorrisMTNS}. However, this multiplier function does not work for all port-Hamiltonian systems with the form \eqref{eq:PHS_1}, as shown by the example of a Timoshenko beam  with parameters from \cite{Mattioni2022}. In the example of a wave equation with constant coefficients, both multiplier functions can be used, but the linear function leads to a better bound on the decay rate. The selection of a multiplier function to optimize the bound on the decay rate is an open research problem.

\bibliographystyle{IEEEtran}
\bibliography{references}

\end{document}